\documentclass{amsart}
\usepackage{comment}
\usepackage{faktor}
\usepackage{array}

\usepackage{tikz-cd}
\usepackage{multicol}

\newcommand{\p}{\textstyle \frac{1}{2}}

\newcommand{\ds}[2]{\delta({#1},{#2})}

\newcommand{\st}[4]
{\ds{#1}{#2} \times \ds{#3}{#4} \rtimes \sigma }

\newcommand{\aaa}{\st{\p}{c}{-a}{b}   }
\newcommand{\bbb}{\st{-a}{b}{\p}{c}   }
\newcommand{\ccc}{\st{-a}{b}{-c}{-\p}   }
\newcommand{\ddd}{\st{-c}{-\p}{-a}{b}   }
\newcommand{\eee}{\st{-c}{-\p}{-b}{a}   }

\newtheorem{subtheorem}{Theorem}[subsection]
\newtheorem{subproposition}[subtheorem]{Proposition}
\newtheorem{sublemma}[subtheorem]{Lemma}

\newtheorem{theorem}{Theorem}[section]
\newtheorem*{theorem*}{Theorem}
\newtheorem{lemma}[theorem]{Lemma}
\newtheorem{corollary}[theorem]{Corollary}

\numberwithin{equation}{section}
\numberwithin{figure}{section}
\numberwithin{table}{section}

\newtheorem{proposition}[theorem]{Proposition}

\theoremstyle{definition}

\theoremstyle{remark}

\usepackage{epsfig,amscd,amssymb,amsxtra,amsmath,amsthm}
\usepackage{url}

\begin{document}

\title[Two linked segments with one half]
{
Induction from two linked segments with  one half border and cuspidal reducibility
}

\author{Igor Ciganovi\'{c}
}
\address{Igor Ciganovi\'{c}, 
Department of Mathematics, Faculty of Science, University of Zagreb, Bijeni\v{c}ka cesta 30, HR-10000 Zagreb, Croatia}
\email{igor.ciganovic@math.hr}

\thanks{This work is 
supported by Croatian Science Foundation under the project
number
HRZZ-IP-2022-10-4615.}
\keywords{Classical group, composition series, induced representations, p-adic field, Jacquet module}
\subjclass[2020]{Primary 22D30, Secondary 22E50, 22D12, 11F85}

\begin{abstract}

 In this paper, we determine the composition series of the induced representation   
 $\delta([\nu^{\frac{1}{2}}\rho,\nu^c\rho])\times 
 \delta([\nu^{-a}\rho,\nu^b\rho])
 \rtimes \sigma$ 
 where 
  $a, b, c \in \mathbb{Z}+\frac{1}{2}$ such that $\frac{1}{2}\leq a < b  < c$,  
 $\rho$ is an irreducible cuspidal unitary representation of a general linear group
and
 $\sigma$ is an irreducible cuspidal representation of a classical group
such that $\nu^\frac{1}{2}\rho\rtimes \sigma$ reduces.
\end{abstract}
\maketitle

\section*{Introduction}

Parabolic induction is an important tool for constructing representations of groups. However, problem of composition 
series of induced representations is solved
only in some special
cases, such as
\cite{bbosnjak:ladder}, 
\cite{degenerate},
\cite{matic-forum},
\cite{muic:composition series}
and
\cite{zelevinsky:ind-repns-II}, and for some low-rank groups.
Motivated by simple results obtained in 
\cite{ciganovic} and its extensions in  
\cite{ciganovic-discrete} and 
\cite{ciganovic-glasnik},
here we continue this effort and calculate composition series of certain induced representations with increased complexity, in the sense that we have more reducibilities occuring when inducing from subsets of starting representations.

To explain this we introduce some notation. 
Let $F$ be a
local non-archimedean field of characteristic different than two, 
$|\ |_F $ its normalized absolute value and $\nu=|\textrm{det}\ |_F $.
Let $\rho$ be
an irreducible cuspidal unitary representation of some $GL(m,F)$, and $x,y\in\mathbb{R}$, such that $y-x+1\in\mathbb{Z}_{> 0}$. 
The set 
$\Delta=[\nu^x \rho,\nu^y \rho]=\{\nu^x \rho,...,\nu^y \rho\}$ is called a segment.
We have a unique irreducible subrepresentation
\[
    \delta=
        \delta(\Delta)=
            \delta([\nu^x \rho,\nu^y \rho])
                \hookrightarrow
                    \nu^y\rho\times \cdots \times \nu^x\rho,
\]
of the parabolically induced representation.
Set 
$e(\delta)=(x+y)/2$.
For a sequence of segments,
such that 
$e(\delta_1)\geq \cdots \geq e(\delta_k)>0$ and
 an irreducible tempered representation $\tau$, of a symplectic or (full) orthogonal group,
 denote by
\[
    L(\delta_1\times\cdots\times \delta_k\rtimes \tau),    
\]
the Langlands quotient of the parabolically induced representation.
  
Now  we fix 
 $\rho$  as above, and so  we shorten the notation 
 $\ds{x}{y}= \delta([\nu^x \rho,\nu^y \rho]$).
Further, let $\sigma$ be an
irreducible cuspidal representation of a symplectic or (full) orthogonal group.
We assume that $\nu^\frac{1}{2}\rho \rtimes \sigma$ reduces.
Let $a, b, c \in \mathbb{Z}+\frac{1}{2}$ such that $\frac{1}{2}\leq a < b  < c$.
In \cite{ciganovic}, we determined composition series of induced representation
\begin{equation*} \label{pr1}
\ds{-b}{c}\times \ds{\p}{a}\rtimes \sigma.
\end{equation*}
We extended this results in \cite{ciganovic-discrete} to certain induced representations of form
\[
     \delta_1 \times \cdots \times \delta_k \rtimes \sigma 
\]
where $\delta_i\times \delta_j$ and 
$\delta_i\times \widetilde{\delta_j}$, $i\neq j$, are irreducible and
$\widetilde{\ }$ stands for the contragredient. 
Loosing this condition on segments, 
in \cite{ciganovic-glasnik},
we determined composition series of 
\[
    \ds{-a}{c} \times \ds{\p}{b} \rtimes \sigma.
\]
In this paper, we determine composition series of 
\[
    \ds{\p}{c} \times \ds{-a}{b} \rtimes \sigma,
\]
having reducibilities 
when
inducing only from segments
in all cases, with or without using contragredient. 
Our main methods are Jacquet modules, M{\oe}glin-Tadić classification of
discrete series and intertwining operators.

To describe the main result of the paper, we introduce some discrete series, appearing as only irreducible subrepresentations in the following induced representations (for more details see Section \ref{sect:2}):
\begin{align*}
&\sigma_a 
    \hookrightarrow 
        \ds{\p}{a}\rtimes \sigma, 
    \textrm{ and similarly for  } 
        \sigma_b \textrm{ and  } \sigma_c,
\\
&\sigma_{b,c}^+ \hookrightarrow \ds{\p}{b}\rtimes\sigma_{c},
    \quad
        \sigma_{b,c}^+ + \sigma_{b,c}^-
            \hookrightarrow 
                \ds{-b}{c}\rtimes \sigma,
\textrm{ and similarly for } \sigma_{a,c}^\pm,
\\
&\sigma_{b,c,a}^\pm
    \hookrightarrow
        \ds{\p}{a} \rtimes \sigma_{b,c}^\pm,
        \quad 
        \sigma_{a,b,c}^+ + \sigma_{a,b,c}^- 
            \hookrightarrow 
                \ds{-a}{b} \rtimes \sigma_c,
\end{align*}
where 
$\sigma_{a,b,c}^+=\sigma_{b,c,a}^+$ denotes the same representation.
Now we have
\begin{theorem*} 
Let $\psi=\ds{\p}{c}\times \ds{-a}{b} \rtimes \sigma$ and define representations
\begin{align*}
    W_1=&\sigma_{b,c,a}^+ +L(\ds{\p}{a}\rtimes \sigma_{b,c}^-),
        \\
    W_2=&L(\ds{\p}{a}\rtimes \sigma_{b,c}^+)
        +
        L(\ds{\p}{b}\rtimes \sigma_{a,c}^-)
        +
        L(\ds{-b}{c}\rtimes \sigma_a)
        +
        L(\ds{-a}{b}\rtimes \sigma_c),
        \\
     W_3=&
            \sigma_{b,c,a}^- +  \sigma_{a,b,c}^-
            +
            L(\ds{\p}{b} \rtimes \sigma_{a,c}^+)
            +
            L(\ds{-a}{c}\rtimes\sigma_{b})+
             L(\ds{-b}{c}\times \ds{\p}{a} 
                                    \rtimes \sigma),
        \\
      W_4=&L(\ds{\p}{b} \times \ds{-a}{c} \rtimes \sigma)
        + L(\ds{\p}{c} \rtimes \sigma_{a,b}^+)
        +L(\ds{\p}{c} \rtimes \sigma_{a,b}^-),
        \\
    W_5=&L(\psi).
\end{align*}
Then there exists a sequence 
$\{0\}=V_0\subseteq V_1 \subseteq V_2 \subseteq V_3
            \subseteq V_4 \subseteq V_5=\psi$,
such that
\begin{equation*}
V_i/V_{i-1}\cong W_i,\quad  i=1,\ldots,5.
\end{equation*}    
\end{theorem*}

The content of the paper is as follows.
After Preliminaries, 
we introduce the notation in Section
\ref{sect:2} and list some reducibility results.
In Section \ref{dekompozicija1} we explain an approach to decompose the induced representation
were an important part plays a kernel of certain intertwining operator.
In Section \ref{O-temperiranim} we provide some results on tempered representations,
that are used in 
Section \ref{diskretni}, when determining discrete
series of the induced representations.
The search for non-tempered candidates is done by 
Section
\ref{netemperirani},
and their multiplicities are determined in
Sections
\ref{confirming-non-tempered subquotients}
and \ref{section-multiplicitet-L--b-c-pola-a},
while Sections 
\ref{faktori-malih-reprezentacija} and
\ref{faktori-velike-reprezentacije}
list composition factors of representations that are considered.
Finally, composition series of needed kernel
are provided in Section
\ref{kompozicijski-niz-malih},
while the main result is proved in Section
\ref{main result}.

\section{Preliminaries}
\label{sect:1}
Our setting is as in \cite{ciganovic-glasnik}, so we briefly recall.
We fix a local non-archimedean field $F$ of characteristic different than two. 
As in \cite{tadic-diskretne}, we fix a tower of symplectic or orthogonal non-degenerate $F$ vector spaces $V_n$, $n\geq 0$ where $n$ is the Witt index. Denote by $G_n$ the group of isometries of $V_n$. It has split rank $n$.  
We fix the set of standard parabolic subgroups
$\{ P_s\}$
in the usual way and have Levi factorization
$ P_s=M_s N_s$.
By Alg $G_n$ we denote smooth representations of $G_n$,
Irr $G_n$ irreducible representations, 
and
subscript $f.l.$ means finite length, $u$ unitary, and $cusp$ cuspidal.
    For $\delta_i \in \textrm{Alg  }GL(n_i,F)$,  $i=1,...,k$ 
and 
$\tau \in \textrm{Alg  }G_{n-m}$ we write
\[
    \delta_1\times\cdots\times\delta_k\rtimes \tau= 
    \text{Ind}_{M_s}^{G_n}( \delta_1\otimes\cdots\otimes\delta_k\otimes \tau)
\]
to denote the normalized parabolic induction.
If $\sigma \in \textrm{Alg }G_n$ we denote by 
$\text{r}_{s}(\sigma)=\text{r}_{M_s}(\sigma)=\text{Jacq}_{M_s}^{G_n}(\sigma)$
the normalized Jacquet module of $\sigma$. 

Let 
$|\ |_F $ be normalized absolute value of $F$ and $\nu=|\textrm{det}\ |_F $.
For
an irreducible cuspidal unitary representation 
$\rho$ of some $GL(m,F)$, and $x,y\in\mathbb{R}$, 
such that 
    $y-x+1\in\mathbb{Z}_{> 0}$, 
the set 
    $\Delta=[\nu^x \rho,\nu^y \rho]=\{\nu^x \rho,...,\nu^y \rho\}$ 
is called a segment.
We have a unique irreducible subrepresentation
\[
    \delta=
        \delta(\Delta)=
            \delta([\nu^x \rho,\nu^y \rho])
                \hookrightarrow
                    \nu^y\rho\times \cdots \times \nu^x\rho,
\]
of the parabolically induced representation.
For $y-x+1\in\mathbb{Z}_{> 0}$ 
define 
$[\nu^x \rho,\nu^y \rho]=\emptyset $
and
$\delta(\emptyset)$ is the irreducible representation 
of the trivial group.
Set 
$e(\delta)=(x+y)/2$.
For a sequence of non-empty segments,
such that 
$e(\delta_1)\geq \cdots \geq e(\delta_k)>0$ and
 an irreducible tempered representation $\tau$, of a symplectic or (full) orthogonal group,
 denote by
\[
    L(\delta_1\times\cdots\times \delta_k\rtimes \tau),    
\]
the Langlands quotient of the parabolically induced representation.

If $\sigma$ is a discrete series representation of $G_n$ then by the
M{\oe}glin-Tadi\'c, now unconditional classification  (\cite{moeglin},\cite{tadic-diskretne}),  
it is described by an admissible triple
\(
(\text{Jord},\sigma_{cusp},\epsilon).
\)

 Let $R(G_n)$ be the free Abelian group generated by classes of irreducible representations of $G_n$.
Put $R(G)=\oplus_{n\geq 0} R(G_n)$.
Let $R^+_0(G)$ be a $\mathbb{Z}_{\geq 0}$
subspan of
classes of irreducible representations.
For $\pi_1,\pi_2\in R(G)$ we define $\pi_1 \leq \pi_2$ if  $\pi_2-\pi_1 \in R^+_0(G)$. 
 Similarly define $R(GL)=\oplus_{n\geq 0}R(GL(n,F))$. We have the map 
 $\mu^* : R(G)\rightarrow R(GL) \otimes R(G)$ defined by
 \[
 \mu^*(\sigma)=1\otimes \sigma + \sum_{k=1}^{n} \text{s.s.}(r_{(k)}(\sigma)),\  \sigma \in R(G_n).
 \]
 where $\text{s.s.}$ denotes the semisimplification.
The following result derives from Theorems 5.4 and 6.5 of  \cite{tadic-structure}, see also section 1.\ in \cite{tadic-diskretne}. They are based on Geometrical  Lemma (2.11 of \cite{bernstein-zelevinsky:ind-repns-I}). 
We use $\widetilde{\ }$ to denote the contragredient.
\begin{theorem} 
Let 
$\sigma$ be a smooth representation of finite length a classical group
and
$[\nu^x\rho,\nu^y\rho]\neq \emptyset$ a segment.
Then
\begin{equation} \label{komnozenje}
\begin{split}
\mu^*(&\delta([\nu^x\rho,\nu^y\rho])\rtimes \sigma)=
\sum_{\delta'\otimes\sigma'\leq \mu^*(\sigma)}
\sum_{i=0}^{y-x+1} \sum_{j=0}^{i}                                \\
&\delta([\nu^{i-y}\widetilde{\rho},\nu^{-x}\widetilde{\rho}])
\times
\delta([\nu^{y+1-j}\rho,\nu^{y}\rho])\times \delta'
\otimes
\delta([\nu^{y+1-i}\rho,\nu^{y-j}\rho])\rtimes \sigma'
\end{split}
\end{equation}
\end{theorem}

Now we write some formulae for Jacquet modules.
Details can be found in \cite{segment} and
 corrections of typographical errors, that we state below, in the Introduction in \cite{segment-ispravljeno}.
Let $\rho$ be an irreducible cuspidal representations of a general linear group, $\sigma$ an irreducible
cuspidal representations of a classical group and 
$c, d \in \mathbb{R}, c + d
    \in \mathbb{Z}_{\geq 0}$.
Assume that 
$\alpha \in \frac{1}{2}\mathbb{Z}_{\geq 0}$ is such that $\nu^\alpha \rho \rtimes\sigma$ reduces. Such $\alpha$ is unique and $\rho$ is selfdual.
Consider induced representation
\[
    \pi=
        \delta([\nu^{-c}\rho,\nu^d \rho])\rtimes \sigma
            \overset{R(G)}{=}
                \delta([\nu^{-d}\rho,\nu^c \rho])\rtimes \sigma.
\]
Three terms are defined: 
$\delta([\nu^{-c}\rho,\nu^d \rho]_+;\sigma)$,
$\delta([\nu^{-c}\rho,\nu^d \rho]_-;\sigma)$
and
\( L_\alpha(\delta([\nu^{-c}\rho,\nu^d \rho]);\sigma)\).
 Each of them is either an irreducible representation or zero.  We have in $R(G)$: 
\begin{equation}
\begin{split}
\label{segment-jacquet-equation}
\delta([\nu^{-a}\rho,\nu^c \rho])\rtimes \sigma
=
&
\delta([\nu^{-a}\rho,\nu^c \rho]_+;\sigma)
+
\delta([\nu^{-a}\rho,\nu^c \rho]_-;\sigma)+
\\
&
L_\alpha(\delta([\nu^{-a}\rho,\nu^c \rho]);\sigma).
\end{split}
\end{equation}
We have
\begin{equation} \label{jacquet-segment}
    \begin{split}
        \mu^*(\delta([\nu^{-c}\rho,\nu^d \rho]_\pm;&\sigma))
        =
        \\
        \sum_{i=-c-1}^{d-1} \sum_{j=i+1}^{d}
        &\delta([\nu^{-i}\rho,\nu^c\rho])
        \times
        \delta([\nu^{j+1}\rho,\nu^d\rho])
        \otimes
        \delta([\nu^{i+1}\rho,\nu^j \rho]_\pm;\sigma)
        +
        \\
        \underset{\textrm{\ } \quad i+j<-1}{
            \sum_{-c-1\leq i\leq c-1}\textrm{\ }
            \sum_{ i+1\leq j\leq c}
            }
        &\delta([\nu^{-i}\rho,\nu^c\rho])
        \times
        \delta([\nu^{j+1}\rho,\nu^d\rho])
        \otimes
        L_\alpha(\delta([\nu^{i+1}\rho,\nu^j \rho]);\sigma)
        \\
        +\sum_{i=-c-1}^{\pm\alpha-1}
        &\delta([\nu^{-i}\rho,\nu^c\rho])
        \times
        \delta([\nu^{i+1}\rho,\nu^d\rho])
        \otimes
        \sigma.
    \end{split}
\end{equation}

The above formula has corrected two typographical errors which exist in \cite{segment}. First, the upper limit in the first sum of the second row needs
to be $d-1$ (instead of $c$, as it is in the published version). Then, the limits of the first
sum in the third row are 
$-c -1\leq  i \leq c-1 $
(instead of $-c -1\leq  i \leq c $ ; the index $c$ does
not give any contribution). 
The same corrections applies to corresponding formulas
in Corollaries 4.3, 5.4 and 6.4 of \cite{segment}.

And for $c<\alpha$ or $\alpha\leq c < d$
we have
\begin{equation} 
\label{jacquet-langlandsov-kvocijent}
    \begin{split}
\mu^*(L_\alpha(\delta([\nu^{-c}\rho,\nu^d \rho&]); \sigma))
    =
        \\
        \underset{\textrm{\ } \quad 0\leq i+j}{
            \sum_{-c-1\leq i\leq d-1}\textrm{\ }
            \sum_{ i+1\leq j\leq d}
            }
         &L(\delta([\nu^{-i}\rho,\nu^c\rho]),
        \delta([\nu^{j+1}\rho,\nu^d\rho]) )
        \otimes
        L_\alpha(\delta([\nu^{i+1}\rho,\nu^j \rho]);\sigma)
        \\
        +\sum_{i=\alpha}^{d} 
        &L(\delta([\nu^{-i}\rho,\nu^c\rho])
        ,
        \delta([\nu^{i+1}\rho,\nu^d\rho]))
        \otimes
        \sigma
    \end{split}
\end{equation}
Also, the above formula has corrected a typographical error existing \cite{segment}: the limits in the first sum in the second row are 
$-c-1 \leq i \leq  d-1$
(instead of
$-c-1 \leq i \leq  d;$  the index $d$ does not contribute in the formula). The same correction applies to corresponding formulas in Corollaries 4.3, 5.4 and 6.4. in \cite{segment}.

In this paper we consider the case $\alpha=\p$. Subquotients of $\pi$ are as follows.

If $ \frac{1}{2} < -c $, then $\pi$ is irreducible and 
$L(\delta([\nu^{-c}\rho,\nu^d \rho]);\sigma)=\pi$.

If $-c \leq \frac{1}{2}  $ then $\pi$ reduces. We denote by 
$\delta([\nu^{-c}\rho,\nu^d \rho]_+;\sigma)$
a unique irreducible subquotient that
has in its minimal standard Jacquet module at least one irreducible subquotient whose all exponents are non-negative.

If $-c=\frac{1}{2}$, then $\pi$ is of length two, 
$\delta([\nu^{-c}\rho,\nu^d \rho]_+;\sigma)$ is a discrete series subrepresentation and
  $L(\delta([\nu^{-c}\rho,\nu^d \rho]);\sigma)$ is the Langlands quotient of  
$\pi$.

If $-c\leq -\frac{1}{2}$, and $c=d$, then $\pi$ is a direct sum of two tempered representations, $\delta([\nu^{-c}\rho,\nu^d \rho]_+;\sigma)$ and $\delta([\nu^{-c}\rho,\nu^d \rho]_-;\sigma)$.

If $-c\leq -\frac{1}{2}$, and $c\neq d $, then $\pi$ is of length three. It has two discrete series representations,
 $\delta([\nu^{-c}\rho,\nu^d \rho]_+;\sigma)$ and $\delta([\nu^{-c}\rho,\nu^d \rho]_-;\sigma)$, and $L(\delta([\nu^{-c}\rho,\nu^d \rho]);\sigma)$ is the Langlands quotient.

\section{Notation and basic reducibilities}
\label{sect:2}
Now we fix the notation and write some reducibility results. 
Let $\rho$ be an irreducible unitary cuspidal representation  of $GL(m_{\rho},F)$ and  
$\sigma$ an irreducible cuspidal representation of $G_n$ such that $\nu^\frac{1}{2}\rho \rtimes \sigma$ reduces. By Proposition 2.4 of
\cite{tadic:red-par-ind} $\rho$ is self-dual. 
We consider 
\[
\frac{1}{2} \leq a, b, c \in \mathbb{Z}+\frac{1}{2},
\]
that need not be fixed, but when appearing together in a formula, we have, depending on which appears, 
$ a < b  < c$. 
We want to decompose the representation 
\[
\psi=\ds{\p}{c} \times \ds{-a}{b} \rtimes \sigma.
\]
First, we shorten some  notations from 
 \eqref{segment-jacquet-equation}:
\begin{align*}
        \sigma_a&=
        \delta([\nu^{\frac{1}{2}}\rho,\nu^a\rho]_+;\sigma),
        &
        \sigma^-_{b,c}&=
        \delta([\nu^{-b}\rho,\nu^c \rho]_-;\sigma),
        &
        \sigma^+_{b,c}&=
        \delta([\nu^{-b}\rho,\nu^c \rho]_+;\sigma),
    \\
        &  
        &
        \tau^+_{a,a}&=
        \delta([\nu^{-a}\rho,\nu^a \rho]_+;\sigma), 
        &
        \tau^-_{a,a}&=
        \delta([\nu^{-a}\rho,\nu^a \rho]_-;\sigma).    
\end{align*}
The following result is Theorem 2.3 from 
\cite{muic:composition series}.
\begin{theorem}
With discrete series being subrepresentations, we have in $R(G)$
\label{muic-diskretne-podreprezentacije} 
\label{prva}
\begin{align*}
\delta([\nu^\frac{1}{2}\rho,\nu^a\rho])\rtimes \sigma
&=
\sigma_a+L(\delta([\nu^\frac{1}{2}\rho,\nu^a\rho])\rtimes \sigma),
\\
\delta([\nu^{-b}\rho,\nu^c\rho])\rtimes \sigma
&=
 \sigma_{b,c}^++
\sigma_{b,c}^-+
L(\delta([\nu^{-b}\rho,\nu^c\rho])\rtimes \sigma).
\end{align*}
\end{theorem}
The next proposition
finishes our
notation. It
follows from Theorem 2.1 of \cite{muic:composition series}.
\begin{proposition}
\label{druga}  
We use $\sigma_{b,c,a}^+=\sigma_{a,b,c}^+$, $\sigma_{b,c,a}^-$ and $\sigma_{a,b,c}^-$
to denote nonisomorphic discrete series subrepresentations,
such that in $R(G)$ we have
\begin{equation*}
\begin{split}
&\delta([\nu^{-b}\rho,\nu^c\rho])\rtimes  \sigma_a=
\sigma_{b,c,a}^+ +
\sigma_{b,c,a}^- +
L(\delta([\nu^{-b}\rho,\nu^c\rho])\rtimes \sigma_a),
\\
&\ds{-a}{b}\rtimes  \sigma_c=
\sigma_{b,c,a}^+ +
\sigma_{a,b,c}^- +
L(\ds{-a}{b}\rtimes  \sigma_c).
\end{split}    
\end{equation*}
\end{proposition}
Next is a consequence of
Theorem 6.3 of \cite{tadic:regular-square}, see also section 3 there.
\begin{proposition} 
\label{multiplicitet-stroge}
We have in $R(G)$, with multiplicity one:
\[
\nu^a\rho\times\cdots \times \nu^{\frac{1}{2}}\rtimes \sigma \geq \sigma_a. 
\]
\end{proposition}

The next lemma follows from Theorem 5.1 of \cite{muic:composition series}, ii)
and Lemma 5.2 of \cite{ciganovic-glasnik}.
\begin{lemma} 
\label{lema-diskretna-podreprezentacija}
    We have in $R(G)$, with discrete series being a 
    subrepresentation
\begin{align} \label{pola-a-sigma-b}
    \ds{\p}{a}\rtimes
      \sigma_{b}
    =&
    \sigma_{a,b}^+
    +
    L(\ds{\p}{a}\rtimes
      \sigma_{b}), 
\\
    \label{pola-b-sigma-a}
    \begin{split}
      \ds{\p}{b}\rtimes
       \sigma_a
     =&
      \sigma_{a,b}^+
     +
      L(\ds{-a}{b}
       \rtimes \sigma)+
      L(
       \ds{\p}{a}\rtimes
        \sigma_b)
      +
      L(
       \ds{\p}{b}\rtimes
        \sigma_a).
    \end{split}    
\end{align}
\end{lemma}
By Proposition 2.4 of  \cite{ciganovic} we have
\begin{lemma} 
\label{lema-pola-a-bc-tvrdnja} 
We have in $R(G)$
    \begin{align} 
    \label{lema-pola-a-bc-f}
    \ds{\p}{a} \rtimes \sigma_{b,c}^\pm
    &=
    \sigma_{b,c,a}^\pm
    +
    L(\ds{\p}{a} \rtimes \sigma_{b,c}^\pm ),
\textrm{ so }
\\
\label{lema-pola-a-bc-f1}
\mu^*(\sigma_{b,c,a}^\pm)
&\geq 
\ds{\p}{a} \otimes \sigma_{b,c}^\pm.
    \end{align}
\end{lemma}
By Propositions 11.6 and 11.8 of \cite{ciganovic-glasnik}
 we have composition series of 
    $\ds{\p}{b}\times \sigma_{a,c}^\pm$.
\begin{lemma} 
\label{faktori-pola-b-sigma -a c pm}
Representations 
$\ds{\p}{b}\times \sigma_{a,c}^+$ 
and
$\ds{\p}{b}\times \sigma_{a,c}^-$
have filtrations
\begin{align*}
    L(\ds{\p}{a}\rtimes \sigma_{b,c}^+)
        +
            L(\ds{-a}{b}\rtimes \sigma_c)
                \hookrightarrow
                    &\ds{\p}{b} \rtimes \sigma_{a,c}^+
                        /
                            \sigma_{a,b,c}^+
                                 \twoheadrightarrow
    L(\ds{\p}{b} \rtimes \sigma_{a,c}^+),
    \\
    L(\ds{\p}{a}\rtimes \sigma_{b,c}^-)
        \hookrightarrow 
            &\ds{\p}{b} \rtimes \sigma_{a,c}^-
                \twoheadrightarrow
                    L(\ds{\p}{b} \rtimes \sigma_{a,c}^-).
\end{align*}
\end{lemma}

By Propositions 3.2 of \cite{faktor} and 
 13.7 of \cite{ciganovic-glasnik} we have
\begin{theorem} 
\label{muic-A3}
Representation $\ds{-a}{c}\rtimes\sigma_{b}$
has a filtration
\begin{align*}
        L(\ds{-b}{c}\rtimes \sigma_a)
        +
        L(\ds{-a}{b}\rtimes \sigma_c)
        \hookrightarrow
        \ds{-a}{c}\rtimes\sigma_{b}/\sigma_{a,b,c}^+
        \twoheadrightarrow
        L(\ds{-a}{c}\rtimes\sigma_{b}).
    \end{align*}
\end{theorem}

\begin{proof}
Composition factors and discrete series being a subrepresentation follow from Propositions 3.2 of \cite{faktor} and its proof.
Proposition 13.7 of \cite{ciganovic-glasnik} relies on 
Proposition 13.1 there, and 
for the sake of completness we write full proof of the latter,
obtaining both results.

We have 
embeddings into a multiplicity one  representation (Theorem 10.3 of \cite{ciganovic-glasnik}) 
\[
    \ds{-a}{c}\rtimes\sigma_{b}
        \hookrightarrow 
            \ds{\p}{b} \times \ds{-a}{c} \rtimes \sigma
                \hookleftarrow
                    \ds{\p}{b} \times \sigma_{a,c}^+
\]
Now 
filtration of $\ds{\p}{b} \times \sigma_{a,c}^+$ 
(Propositions 11.6 of \cite{ciganovic-glasnik}) 
show that 
$L(\ds{-a}{b}\rtimes \sigma_c)$ 
needs quotient
$\ds{-a}{c}\rtimes\sigma_{b}/\sigma_{b,c,a}^+$ to embed.

Next, consider
following two compositions of an embedding and an epimorphism 
\begin{align*}
    \ds{-a}{c}\rtimes\sigma_{b}
        \hookrightarrow 
            \ds{\p}{c} \times \ds{-a}{-\p} \rtimes \sigma_{b} 
                \twoheadrightarrow
                    \ds{\p}{c} \times \sigma_{a,b}^+,
\\
    \ds{-b}{c} \times \sigma_a/ \sigma_{b,c,a}^-
        \hookrightarrow 
           \ds{\p}{c} \times \ds{-b}{-\p} \times \sigma_a/ \sigma_{b,c,a}^-
            \twoheadrightarrow
                \ds{\p}{c} \times \sigma_{a,b}^+.
\end{align*}
By Lemma 9.1 of \cite{ciganovic-glasnik} 
and 
\eqref{lema-pola-a-bc-f1}
we have 
$\ds{\p}{a}\otimes L(\ds{-b}{c}\rtimes \sigma) 
\leq \mu^*(L(\ds{-b}{c}\rtimes \sigma_a))$
and 
$\ds{\p}{a} \otimes \sigma^+_{b,c}\leq \mu^*(\sigma^+_{b,c,a})$
and one can check that they both 
appear with multiplicity one in all representations in 
above two compositions.
Thus both compositions are embeddings and they 
induce an embedding
\[
    \ds{-b}{c} \times \sigma_a/ \sigma_{b,c,a}^- 
        \hookrightarrow
            \ds{-a}{c}\rtimes\sigma_{b}.
\]
This shows that  
$ L(\ds{-b}{c}\rtimes \sigma_a)$
needs quotient
$\ds{-a}{c}\rtimes\sigma_{b}/\sigma_{b,c,a}^+$ to embed.
\end{proof}

The last two statements are Proposition 12.1 of \cite{ciganovic-glasnik} and the main result there.
\begin{theorem} \label{ciganovic-A2} 
There exists a filtration $\{ V_i\}$ of 
$
\ds{-c}{b} \times 
\ds{\p}{a} \rtimes \sigma
$
such that
\begin{align*}
V_1 &\cong     L(\ds{-b}{c}\rtimes \sigma_a),
\\
V_2/V_1 & \cong
    \sigma_{b,c,a}^+
    +
    \sigma_{b,c,a}^-
    +
    L(\ds{-b}{c} \times \ds{\p}{a} \rtimes \sigma),
\\
V_3/V_2 & \cong 
    L(\ds{\p}{a}\rtimes \sigma_{b,c}^+)
    +
    L(\ds{\p}{a}\rtimes \sigma_{b,c}^-).
\end{align*}
\end{theorem}

\begin{theorem}
\label{faktori-pola-b-ds-a-c}
There exists a filtration 
$\{V_i\}$
of
 $\ds{-a}{c}\times \ds{\p}{b}\rtimes \sigma$ 
 such that
\begin{align*}
V_1 \cong &
    \sigma_{b,c,a}^+ 
    +L(\ds{\p}{a}\rtimes \sigma_{b,c}^-),
        \\
V_2/V_1 \cong &
    L(\ds{\p}{a}\rtimes \sigma_{b,c}^+)
    +
    L(\ds{-a}{b}\rtimes \sigma_c)
    +
        L(\ds{\p}{b}\rtimes \sigma_{a,c}^-)
        +
        L(\ds{-b}{c}\rtimes \sigma_a),
        \\
V_3/V_2\cong & 
    L(\ds{\p}{b} \rtimes \sigma_{a,c}^+)+
            L(\ds{-a}{c}\rtimes\sigma_{b})+
        \sigma_{b,c,a}^- + 
                            L(\ds{-b}{c}\times
                                \ds{\p}{a} 
                                    \rtimes \sigma),
        \\
V_4/V_3\cong  &
    L(\ds{-a}{c}\times \ds{\p}{b}\rtimes \sigma).    
\end{align*}

\end{theorem}


\section{An estimate for 
$\ds{\p}{c}\times \ds{-a}{b}\rtimes \sigma$
} 
\label{dekompozicija1}
We start our decomposition by considering some standard intertwining operators.
\begin{alignat*}{3}
\aaa 
    &\overset{f_1}{\longrightarrow}
\bbb
    &&\overset{f_2}{\longrightarrow}
\ccc
\\    
    &\overset{f_3}{\longrightarrow}
\ddd
    &&\overset{f_4}{\longrightarrow}
\eee.
\end{alignat*} 
For all $i\geq 1$ denote
$K_i=\textrm{Ker }f_i$.
By Theorem  \ref{prva} 
we have 
\begin{align*}
&K_1\cong      \ds{\p}{b}\times \ds{-a}{c}\rtimes \sigma,
&&
K_2\cong \ds{-a}{b}\rtimes \sigma_c,
 \\
 &K_3\cong
 \ds{-c}{b} \times \ds{\p}{a}\rtimes \sigma,
 &&
 K_4 \cong
\ds{-c}{-\p}\rtimes\sigma_{a,b}^+
      +
      \ds{-c}{-\p}\rtimes\sigma_{a,b}^-. 
\end{align*}
We denoted $\psi=\ds{\p}{c}\times \ds{-a}{b}\rtimes \sigma$ and obtain the estimate, in $R(G)$:
\begin{equation}
    \label{long-intertwining}
\forall i\quad   K_i \leq  \psi \leq K_1+K_2+K_3 +K_4 +L(\psi).
\end{equation}
By Section \ref{sect:2}, 
it remains to decompose $K_4$ to obtain all composition factors of $\psi$.


\section{On some irreducible tempered repesentations}
\label{O-temperiranim}
Here we provide some results that we use for multiplicities of 
discrete series and their Jacquet modules. 
By \cite{goldberg-tempered},  
we have decomposition of the induced representation 
\begin{align}
\label{temperiraneT1}
    \ds{-a}{a} \rtimes \sigma_c &= T^+ + T^-
\end{align}
into a direct sum of non-equivalent tempered 
representations.
Similar to Definition 4.6 of \cite{tadic:tempered},
we use $T^+$ to denote the unique irreducible subquotient 
such that
\begin{align*}
    \mu^*(T^+)& \geq \ds{\p}{a}  \times \ds{\p}{a} \otimes \pi,
\end{align*}
for some irreducible representation $\pi$.
Use \eqref{komnozenje} and 
\eqref{jacquet-segment}
to calculate
\begin{align}
\label{temperiraneT2}
    \mu^*(T^+)& \geq \ds{\p}{a}  \times \ds{\p}{a} \otimes \sigma_c.
\end{align}
The following lemma gives some more information. 
\begin{lemma} 
\label{jacquetovi-moduli-od-T}
Writting out all irreducible subquotients of form 
$\ds{\p}{a} \otimes \ldots$,
with maximum multiplicities, we have
\begin{align*}
    \mu^*(T^+)& \geq
        \ds{\p}{a} \otimes L(\ds{\p}{a} \rtimes \sigma_c)
        +
        2 \cdot \ds{\p}{a}  \otimes \sigma_{a,c}^+, 
    \\
    \mu^*(T^-)& \geq
        \ds{\p}{a} \otimes L(\ds{\p}{a} \rtimes \sigma_c).
\end{align*}

\end{lemma}

\begin{proof}
    Looking for $\ds{\p}{a} \otimes \ldots $ in \eqref{temperiraneT1} we obtain
    \[
          \mu^*(\ds{-a}{a} \rtimes \sigma_c)  
            \geq    
                2 \cdot \ds{\p}{a} 
                    \otimes 
                        \ds{\p}{a} \rtimes \sigma_c,
    \]
    where
    by \eqref{pola-a-sigma-b}, in $R(G)$: 
    $\ds{\p}{a} \rtimes \sigma_c
        =\sigma_{a,c}^+
            +
                L(\ds{\p}{a} \rtimes \sigma_c)$.
    Since 
    $\mu^*(\sigma_{a,c}^+)\geq \ds{\p}{a}\otimes \sigma_c$, transitivity of Jacquet module and \eqref{temperiraneT2} imply
    $\mu^*(T^+) \geq
        2 \cdot \ds{\p}{a}  \otimes \sigma_{a,c}^+$.
        The rest of the claim follows as
        \[
         T^+ + T^-=
         \ds{-a}{a} \rtimes \sigma_c 
         \hookrightarrow
         \ds{\p}{a}
         \times \ds{-a}{-\p} \times \sigma_c
        \]
        implies 
        $\mu^*(T^\pm) \geq 
            \ds{\p}{a} \otimes 
                \ds{-a}{-\p} \otimes \sigma_c$
        and $\mu^*(\sigma_{a,c}^+)\ngeq \ds{-a}{-\p} \otimes \sigma_c$. 
\end{proof}

\begin{proposition}
    \label{temperirane-u-pola-c-a-a}
    Writting with maximum multiplicities we have in $R(G)$:
    \begin{align}
        \label{temperirane-u-pola-c-a-a-1}
        \ds{\p}{c}\times \ds{-a}{a} \rtimes \sigma
        \geq 1 \cdot T^+ +  1 \cdot T^-,
    \end{align}
    where $\ds{-a}{a} \rtimes \sigma_c = T^+ + T^-$ and 
    $\mu^*(T^+) \geq \ds{\p}{a}  \times \ds{\p}{a} \otimes \sigma_c$.
\end{proposition}

\begin{proof}
    Since $\ds{-a}{a} \rtimes \sigma_c \leq 
           \ds{\p}{c}\times \ds{-a}{a} \rtimes \sigma $, 
    we look for $\ds{\p}{a}\otimes \ldots $ in 
    $\mu^*(\ds{\p}{c}\times \ds{-a}{a} \rtimes \sigma)$ and obtain
    \[
        2\cdot \ds{\p}{a}\otimes 
            \ds{\p}{a} \times \ds{\p}{c} \rtimes \sigma. 
    \]
    By Lemma 8.3 of \cite{ciganovic-glasnik},
    $L(\ds{\p}{a} \rtimes \sigma_c)$ appears once in 
    $\ds{\p}{a} \times \ds{\p}{c} \rtimes \sigma$, so 
    Lemma
    \ref{jacquetovi-moduli-od-T}  gives 
    \eqref{temperirane-u-pola-c-a-a-1}.
\end{proof}

\section{Discrete series subquotients}
\label{diskretni}
Here we determine discrete series subquotients in three induced representations
\begin{equation} \label{A-A1}
\begin{split}
\ds{\p}{c} \times \ds{-a}{b}\rtimes \sigma 
\geq 
      \ds{\p}{c} \rtimes \sigma_{a,b}^+
      +
      \ds{\p}{c}\rtimes \sigma_{a,b}^-.
\end{split}
\end{equation}
Note that by Lemma 13.4 of 
\cite{ciganovic-glasnik} we know 
$\sigma_{b,c,a}^+ 
\hookrightarrow \ds{\p}{c} \rtimes \sigma_{a,b}^+$.
\begin{proposition}
\label{svi-diskretni-u-pola-c--a-b-sigma}
    Writting all discrete series 
    with maximum multiplicities 
    we have 
    \begin{align}
        \label{diskretne-u-pola-c-a-b1}
        \ds{\p}{c} \times \ds{-a}{b}\rtimes \sigma 
            &\geq 1 \cdot \sigma_{b,c,a}^+ 
                    + 1 \cdot \sigma_{b,c,a}^- 
                    + 1 \cdot \sigma_{a,b,c}^-,
        \\
        \label{diskretne-u-pola-c-a-b2}
        \ds{\p}{c} \rtimes \sigma_{a,b}^+
            &\geq 
                1 \cdot \sigma_{b,c,a}^+ 
                    + 1 \cdot \sigma_{b,c,a}^-,
        \\
        \label{diskretne-u-pola-c-a-b3}
        \ds{\p}{c}\rtimes \sigma_{a,b}^-
            &\geq  
                1 \cdot \sigma_{a,b,c}^-.
    \end{align}
\end{proposition}
\begin{proof}
    As in Lemma 5.1 of \cite{ciganovic-glasnik}, we obtain that 
    $\sigma_{b,c,a}^+=\sigma_{a,b,c,}^+$, $\sigma_{b,c,a}^-$ and
    $\sigma_{a,b,c}^-$ are only possible discrete subquotients in all equations. 
    Proposition \ref{druga} and Section \ref{dekompozicija1} imply that they do appear in 
    \eqref{diskretne-u-pola-c-a-b1}.
    
    To check multiplicity one of 
    $\sigma_{b,c,a}^\pm$ in 
    \eqref{diskretne-u-pola-c-a-b1}, use 
     Lemma 5.2 of \cite{ciganovic-glasnik} and 
     \eqref{lema-pola-a-bc-f1} 
      to see multiplicity one  of
    $\ds{\p}{a} \otimes \sigma_{b,c}^\pm$ in the appropriate Jacquet module. 
    Now, multiplicity  of $\sigma_{b,c,a}^\pm$
    is one in \eqref{diskretne-u-pola-c-a-b2}, by
    Lemma \ref{lema-diskretna-podreprezentacija}, since
    $\mu^*(\sigma_{a,b}^+)
        \geq \ds{\p}{a} \otimes \sigma_{b}$. 
    Thus multiplicity  of $\sigma_{b,c,a}^\pm$ is zero in 
        \eqref{diskretne-u-pola-c-a-b3}.

    To check multiplicity  of $\sigma_{a,b,c}^\pm$, we
    search for 
    
    $\ds{a+1}{b}\otimes\ldots$ in 
    $\mu^*(\ds{\p}{c} \times \ds{-a}{b}\rtimes \sigma)$ and obtain
    \begin{align*}
        \ds{a+1}{b}\otimes 
            \ds{\p}{c}\times 
                \ds{-a}{a}\rtimes
                    \sigma.
    \end{align*}
    Proposition \ref{temperirane-u-pola-c-a-a}
    implies multiplicity one of $\sigma_{a,b,c}^\pm$ in  
    \eqref{diskretne-u-pola-c-a-b1}. Here we had contribution
    \begin{align*}
    \mu^*(\ds{-a}{b}\rtimes \sigma)
        &\geq
            \ds{a+1}{b}\otimes \ds{-a}{a}\rtimes\sigma
        \\    
        &=
            \ds{a+1}{b}\otimes \tau_{a,a}^+
                +
                    \ds{a+1}{b}\otimes \tau_{a,a}^-.
    \end{align*}
     Since
    $\mu^*(\sigma_{a,b}^\pm)
        \geq
            \ds{a+1}{b} \otimes \tau_{a,a}^\pm$,
    the contribution comes from 
        $\sigma_{a,b}^+ +\sigma_{a,b}^-$.
    Using
    $\ds{\p}{c} \rtimes \sigma_{a,b}^+
        \geq 
            \sigma_{b,c,a}^+$,       
    we have
    $\ds{\p}{c} \rtimes \sigma_{a,b}^\pm
        \geq 
            \sigma_{a,b,c}^\pm$,
    proving \eqref{diskretne-u-pola-c-a-b2} and \eqref{diskretne-u-pola-c-a-b3}.
\end{proof}


\section{Non-tempered candidates}
\label{netemperirani}
As noted in Section \ref{dekompozicija1}, we  search for non-tempered subquotients
in 
\[
	\ds{\p}{c}\rtimes
     	 \sigma_{a,b}^+
      				\quad\textrm{ and }\quad
      \ds{\p}{c}\rtimes
     	 \sigma_{a,b}^-.
\]

\begin{lemma} \label{veza-diskretnih}
With maximum multiplicity of $\sigma_{b,c}^\pm$ being one, and $\epsilon,\eta \in \{\pm \}
$ 
we have
  \begin{equation*}
     \begin{split}
     \sigma_{b,c}^\epsilon \leq \ds{a+1}{c}\rtimes 
			\sigma_{a,b}^\eta\quad  \iff \epsilon=\eta,
      \\
      \sigma_{b,c}^+ + \sigma_{b,c}^- \leq 
	\ds{a+1}{c} \times \ds{-b}{a}\rtimes \sigma.
    \end{split}
  \end{equation*}  
\end{lemma}

\begin{proof}
Denote the last representation by $\pi$. We have embeddings
		\begin{align*}
			\qquad \qquad \qquad
			\sigma_{b,c}^\epsilon
			\rightarrow
			\ds{-b}{c}\rtimes \sigma  
			\rightarrow \pi
		\end{align*}
By \eqref{jacquet-segment}, $\ds{a+1}{c}\times\ds{a+1}{b}\otimes\sigma_{a,a}^\epsilon$
appears once in $\mu^*(\sigma_{b,c}^\epsilon)$.
Check 
the same for
$\mu^*(\pi)$. 
The first claim follows as
$\mu^*(\sigma_{a,b}^\epsilon)\geq \ds{a+1}{b}\otimes \sigma_{a,a}^\epsilon$.
\end{proof}

Now we determine non-tempered candidates.
\begin{proposition} 
\label{nediskretni-kandidati-u-pozitivnoj-C}
Fix $\epsilon=+$ or $-$.
If $\pi$ is a non-tempered subquotient of  
$\ds{\p}{c}\rtimes\sigma_{a,b}^\epsilon$, different from its Langlands quotient, 
		then $\pi$ can be 
	\[
			L(\ds{\p}{a}\rtimes \sigma_{b,c}^\epsilon) 
		\quad  \textrm{   or } \quad
    		  	L(\ds{\p}{b}\rtimes \sigma_{a,c}^\epsilon),
	\]
	and moreover if $\epsilon=+$, then $\pi$ can also be
	\[
			L(\ds{-a}{c}\rtimes \sigma_b), \quad 
			 L(\ds{-a}{b}\rtimes \sigma_c)
		\quad  \textrm{   or } \quad
			L(\ds{-b}{c}\rtimes\sigma_a).
	\]
\end{proposition}

\begin{proof}
We use Lemma 2.2 of \cite{muic:composition series}
(in terms of that lemma 
$\pi \leq \ds{-l_1}{l_2}\rtimes \sigma$, 
$-l_1=\frac{1}{2}$, $l_2=c$ and $\sigma=\sigma_{a,b}^\epsilon$).
So we look for possible embeddings
	\begin{equation} \label{moguca_ulaganja}
		\pi \hookrightarrow \ds{-\alpha_1}{\beta_1}\rtimes \pi',
	\end{equation}
where $-\alpha_1+\beta_1<0$ and $\pi'$ is irreducible.
By the lemma, there exists an irreducible representation $\sigma_1$ such that
	\begin{equation} \label{moguca_ulaganja_uvjet_mu}
    		\begin{cases}
    			\mu^*(\sigma_{a,b}^\epsilon)\geq \ds{\p}{\beta_1}\otimes \sigma_1
    			\\
    			\pi'\leq \ds{\alpha_1+1}{c}\rtimes \sigma_1
    		\end{cases}
	\end{equation}
and we must have
	\begin{equation} \label{moguca_ulaganja_uvjet_nejednakost}
		\begin{cases}
			-\frac{1}{2} \leq \beta_1  \\
			c \geq  \alpha_1  > \beta_1, -\frac{1}{2}  \\
			\alpha_1 \geq \p.
		\end{cases}
	\end{equation}
We have two possible cases:
\begin{itemize}
    \item[a)]
    $\beta_1=-\frac{1}{2}$. Now $\sigma_1=\sigma_{a,b}^\epsilon$.
    \begin{itemize}
    		\item[$\bullet$] Assume that $\pi'$ is tempered. 
		We may take $2\alpha_1+1\in \textrm{Jord}_\rho(\sigma_{a,b}^\epsilon)$. 				So
	\begin{itemize}
		\item[$1)$] $\alpha_1=a$. Then		
		\[
    			\pi'\leq \ds{a+1}{c}\rtimes \sigma_{a,b}^\epsilon.
    		\]
		 Lemma 8.1 of 
		\cite{tadic-diskretne} implies  
		$\pi'=\sigma_{b,c}^+$ or $\sigma_{b,c}^-$. 
		By Lemma \ref{veza-diskretnih} we have $\pi'\cong \sigma_{b,c}^\epsilon$.  
		Thus $\pi\cong L(\ds{\p}{a}\rtimes \sigma_{b,c}^\epsilon)$.
		\item[$2)$] $\alpha_1=b$. 
			Then		
    			$\pi'\leq \ds{b+1}{c}\rtimes \sigma_{a,b}^\epsilon$.
			 Lemma 8.1 of 
		\cite{tadic-diskretne} implies  
		$\pi'=\sigma_{a,c}^+$ or $\sigma_{a,c}^-$. 
		By Section 8. of  \cite{tadic:tempered}, we have $\pi'=\sigma_{a,c}^\epsilon$. So  
		$\pi \cong L(\ds{\p}{b}\rtimes \sigma_{a,c}^\epsilon)$.
	\end{itemize}    		

    \item[$\bullet$] If $\pi'$ is not tempered, by
    Lemma 2.2 of \cite{muic:composition series}, 
	there exist
 	$\alpha_2\in \mathbb{Z} +\p$,
    	$2\beta_2+1 \in \textrm{Jord}_\rho(\sigma_{a,b}^\epsilon)
	$ 
	and 
	$(2\beta_2+1)_-=2(\beta_2)_-+1
 	\in \textrm{Jord}_\rho(\sigma_{a,b}^\epsilon)=\{2a+1,2b+1 \}$,
    such that 
    (in terms of the lemma $\alpha_1\leq (\beta_2)_-<\beta_2<\alpha_2\leq l_2$
		and 	$-\alpha_1+\beta_1\leq -\alpha_2+\beta_2$)    we 	have
    	\[
		\begin{cases}
   			 \alpha_1 \leq a <  b <\alpha_2 \leq c  \\
			-\alpha_1 -\p   \leq  -\alpha_2 -b.
    		\end{cases}
	\]
		The second equation gives $\alpha_2-\alpha_1\leq   - (\p -b) <0$, 
				which contradicts the first one.
    \end{itemize}

    \item[b)] $\beta_1 >-\frac{1}{2}$. Then, by the lemma, 
    $2\beta_1+1\in \textrm{Jord}_\rho(\sigma_{a,b}^\epsilon)=\{2a+1,2b+1 \}$.
    Here we have two options:
	\begin{itemize}
		\item[$\bullet$] 	$\beta_1=a$. 
							Then 
					$\mu^*(\sigma_{a,b}^\epsilon)\geq \ds{\p}{a}\otimes \sigma_1$.
								Now \eqref{jacquet-segment} implies $\epsilon=+$,  
					$\sigma_1=\sigma_b$ and $\pi'\leq \ds{\alpha_1+1}{c}\rtimes \sigma_b$
							for some $a<\alpha_1 \leq c$; $\alpha_1\in \mathbb{Z}+\p$.
		By Proposition 3.1 of \cite{muic:composition series}, depending on $\alpha_1$, $\pi'$ can be:
\begin{equation*}\quad \qquad
			\pi'\cong
		\left \{
\begin{alignedat}{2}
L(\ds{\alpha_1+1}{b}\rtimes \sigma_c) \textrm{ or }  L(\ds{\alpha_1+1}{c}\rtimes \sigma_b),  					 && a <\alpha_1 & <b, \\
\sigma_c  \textrm{ or } L(\ds{b+1}{c}\rtimes \sigma_b), 
				 && \alpha_1 & =b, \\
L(\ds{\alpha_1+1}{c}\rtimes \sigma_b), 
				 && \quad b <\alpha_1 & <c, \\
\sigma_b,  && \alpha_1 & =b. 
\end{alignedat}
\right.
\end{equation*}
		Assume any of cases $\pi'\cong L(\ds{\alpha_1+1}{c}\rtimes \sigma_b)$, where $a<\alpha_1<c$.
		Introducing $\Phi$, we have embeddings
		\begin{align*}
				&\pi \hookrightarrow \ds{-\alpha_1}{a} 
						\times \ds{-c}{-\alpha_1-1} \rtimes \sigma_b =: \Phi,
		\\
		 		&L(\ds{-a}{c}\rtimes \sigma_b ) \hookrightarrow
					\ds{-c}{a}\rtimes \sigma_b \hookrightarrow \Phi.
		\end{align*}
		By 9.1 of \cite{zelevinsky:ind-repns-II}, and transitivity of Jacquet modules, for every
		\[
				\ds{-\alpha_1}{a} 
						\otimes \ds{-c}{-\alpha_1-1} \otimes \sigma_b \leq \mu^*(\Phi) 
		\]		
		there also exists one 
					\[
							\ds{-c}{a} \otimes \sigma_b \leq \mu^*(\Phi).
					\]
		It is not hard to check that the last multiplicity is one. 

		Thus $\pi \cong L(\ds{-a}{c}\rtimes \sigma_b )$, as in the case $\pi'\cong\sigma_b$.
		Other cases similarly give $\pi \cong L(\ds{-a}{b}\rtimes \sigma_c )$.

		\item[$\bullet$] $\beta_1=b$.  
							Then 
							$\mu^*(\sigma_{a,b}^\epsilon)\geq \ds{\p}{b}\otimes \sigma_1$.
								Now \eqref{jacquet-segment} implies $\epsilon=+$,  
							 $\sigma_1=\sigma_a$ and 
						$\pi'\leq \ds{\alpha_1+1}{c}\rtimes \sigma_a$
							for some $b<\alpha_1 \leq c$, $\alpha_1\in \mathbb{Z}+\p$,
							which is irreducible by Proposition 3.1 of \cite{muic:composition series}.
			Thus	
			\[
			\pi\hookrightarrow \ds{-\alpha_1}{b} \times \ds{-c}{-\alpha_1-1}\rtimes \sigma_a.
			\]
			Similarly as above, we obtain $\pi \cong L(\ds{-b}{c}\rtimes \sigma_a)$.
	\end{itemize}
	\end{itemize}
\end{proof}


\section{Non-tempered subquotients and their multiplicities}
\label{confirming-non-tempered subquotients}

Now we show multiplicity one for non-tempered candidates 
from
Proposition 
    \ref{nediskretni-kandidati-u-pozitivnoj-C}.

\subsection{
Multiplicity of
$L(\ds{\p}{a}\rtimes \sigma_{b,c}^\pm)$
}
\begin{sublemma}
    Writting with maximum multiplicity,
    we have in $R(G)$:
    \begin{align*} 
        \mu^*(\ds{\p}{c} \times \ds{-a}{b}\rtimes \sigma)
        \geq
            1 \cdot \ds{-a}{-\p}\otimes \sigma_{b,c}^\pm,
        \\
        \mu^*(\ds{\p}{c}\rtimes \sigma_{a,b}^\pm)
        \geq
            1 \cdot \ds{-a}{-\p}\otimes \sigma_{b,c}^\pm.
    \end{align*}
\end{sublemma}

\begin{proof}
    Looking for $\ds{-a}{-\p}\otimes \ldots$
     in $\mu^*(\ds{\p}{c} \times \ds{-a}{b}\rtimes \sigma)$
     we obtain 
     \begin{align*}
        \ds{-a}{-\p}\otimes 
         \ds{a+1}{c} \times \ds{-a}{b}\rtimes \sigma,
     \end{align*}
     with contribution
     $\mu^*(\ds{-a}{b}\rtimes \sigma)\geq 1 \otimes \ds{-a}{b}\rtimes \sigma$.
     We prove both claims using Lemma \ref{veza-diskretnih}.
\end{proof}

\begin{subproposition}
\label{mult-L1}
Writting with maximum multiplicity,
    we have in $R(G)$:
    \begin{align*} 
        \ds{\p}{c} \times \ds{-a}{b}\rtimes \sigma
        \geq
            1 \cdot L( \ds{\p}{a}\rtimes \sigma_{b,c}^\pm),
        \\
        \ds{\p}{c}\rtimes \sigma_{a,b}^\pm
        \geq
           1 \cdot L( \ds{\p}{a}\rtimes \sigma_{b,c}^\pm).
    \end{align*}
\end{subproposition}


\subsection{
Multiplicity of
$L(\ds{\p}{b}\rtimes \sigma_{a,c}^\pm)$
}

\begin{sublemma}
\label{diskretne-u-b-i-pola-a-c}
  Writting with maximum multiplicity,
    we have in $R(G)$:
  \begin{equation*}
     \begin{split}
1 \cdot     \sigma_{b,c}^\pm \leq \ds{b+1}{c}\rtimes 
\sigma_{a,b}^\pm,
      \\
1 \cdot      \sigma_{b,c}^+ + 1 \cdot \sigma_{b,c}^- \leq \
        \ds{b+1}{c}
            \times \ds{-a}{b}\rtimes \sigma.
    \end{split}
  \end{equation*}  
\end{sublemma}
\begin{proof}
This follows from  
Section 8.\ of  \cite{tadic:tempered} and  \eqref{komnozenje}.
\end{proof}

\begin{sublemma}
    Writting with maximum multiplicity,
    we have in $R(G)$:
    \begin{align*} 
        \mu^*(\ds{\p}{c} \times \ds{-a}{b}\rtimes \sigma)
        \geq
            1 \cdot \ds{-b}{-\p}\otimes \sigma_{a,c}^\pm,
        \\
        \mu^*(\ds{\p}{c}\rtimes \sigma_{a,b}^\pm)
        \geq
            1 \cdot \ds{-b}{-\p}\otimes \sigma_{a,c}^\pm.
    \end{align*}
\end{sublemma}

\begin{proof}
    Looking for $\ds{-b}{-\p}\otimes \ldots$
     in $\mu^*(\ds{\p}{c} \times \ds{-a}{b}\rtimes \sigma)$
     we obtain 
     \begin{align*}
        \ds{-b}{-\p}\otimes 
         \ds{b+1}{c} \times \ds{-a}{b}\rtimes \sigma,
     \end{align*}
     with contribution
     $\mu^*(\ds{-a}{b}\rtimes \sigma)\geq 1 \otimes \ds{-a}{b}\rtimes \sigma$.
     We prove both claims using 
     Lemma \ref{diskretne-u-b-i-pola-a-c}.
\end{proof}

\begin{subproposition}
\label{mult-L2}
Writting with maximum multiplicity,
    we have in $R(G)$:
    \begin{align*} 
        \ds{\p}{c} \times \ds{-a}{b}\rtimes \sigma
        \geq
            1 \cdot L( \ds{\p}{b}\rtimes \sigma_{a,c}^\pm),
        \\
        \ds{\p}{c}\rtimes \sigma_{a,b}^\pm
        \geq
           1 \cdot L( \ds{\p}{b}\rtimes \sigma_{a,c}^\pm).
    \end{align*}
\end{subproposition}

\subsection{
Multiplicity of
$L(\ds{-a}{c}\rtimes \sigma_{b})$
}

\begin{sublemma}
    Writting with maximum multiplicity,
    we have in $R(G)$:
    \begin{align*} 
        \mu^*(\ds{\p}{c} \times \ds{-a}{b}\rtimes \sigma)
        \geq
            1 \cdot \ds{-c}{a}\otimes \sigma_{b},
        \\
        \mu^*(\ds{\p}{c}\rtimes \sigma_{a,b}^+)
        \geq
            1 \cdot \ds{-c}{a}\otimes \sigma_{b}.
    \end{align*}    
\end{sublemma}
\begin{proof}
Looking for $\ds{-c}{a}\otimes \ldots$
     in $\mu^*(\ds{\p}{c} \times \ds{-a}{b}\rtimes \sigma)$
     we obtain 
     \begin{align*}
        \ds{-c}{-\p} \times  \ds{\p}{a} 
        \otimes 
         \ds{\p}{b}\rtimes \sigma,
     \end{align*}
     with contribution
     $\mu^*(\ds{-a}{b}\rtimes \sigma)\geq 
     \ds{\p}{a}
     \otimes \ds{\p}{b}\rtimes \sigma$.
    The first claims followby Theorem  
    \ref{muic-diskretne-podreprezentacije}
    and the second by 
        $\mu^*(\sigma_{a,b}^+)\geq 
            \ds{\p}{a}\otimes \sigma_{b}$.
\end{proof}

\begin{subproposition}
\label{mult-L3}
Writting with maximum multiplicity,
    we have in $R(G)$:
    \begin{align*} 
        \ds{\p}{c} \times \ds{-a}{b}\rtimes \sigma
        \geq
            1 \cdot L(\ds{-a}{c}\rtimes \sigma_{b}),
        \\
        \ds{\p}{c}\rtimes \sigma_{a,b}^+
        \geq
           1 \cdot L(\ds{-a}{c}\rtimes \sigma_{b}).
    \end{align*}
\end{subproposition}

\subsection{Multiplicities of
$L(\ds{-a}{b}\rtimes \sigma_c)$ 
and 
$L(\ds{-b}{c}\rtimes \sigma_a)$
}

\begin{sublemma}
    Writting with maximum multiplicity,
    we have in $R(G)$:
    \begin{align*} 
        \mu^*(\ds{\p}{c} \times \ds{-a}{b}\rtimes \sigma)
        \geq
            1 \cdot \ds{\p}{a}\otimes 
                    L(\ds{\p}{b}\rtimes \sigma_c)
            +
            1 \cdot \ds{\p}{a}\otimes 
                    L(\ds{-b}{c}\rtimes \sigma).
        \\
        \mu^*(\ds{\p}{c}\rtimes \sigma_{a,b}^+)
        \geq
            1 \cdot \ds{\p}{a}\otimes 
                    L(\ds{\p}{b}\rtimes \sigma_c)
            +
            1 \cdot \ds{\p}{a}\otimes 
                    L(\ds{-b}{c}\rtimes \sigma).
    \end{align*}
\end{sublemma}
\begin{proof}
    Looking for $\ds{\p}{a}\otimes \ldots$
     in $\mu^*(\ds{\p}{c} \times \ds{-a}{b}\rtimes \sigma)$
     we obtain 
     \begin{align*}
        \ds{\p}{a}\otimes 
         \ds{\p}{c} \times \ds{\p}{b}\rtimes \sigma,
     \end{align*}
     with contribution
     $\mu^*(\ds{-a}{b}\rtimes \sigma)
        \geq 
            \ds{\p}{a} \otimes \ds{\p}{b}
                \rtimes \sigma$.
        The first claims follows by Lemmas 5.2 and 8.3 of 
        \cite{ciganovic-glasnik} and the second by 
        $\mu^*(\sigma_{a,b}^+)\geq 
            \ds{\p}{a}\otimes \sigma_{b}$.
\end{proof}

\begin{subproposition}
\label{mult-L4}
Writting with maximum multiplicity,
    we have in $R(G)$:
    \begin{align*} 
        \ds{\p}{c} \times \ds{-a}{b}\rtimes \sigma
        &\geq
            1 \cdot L(\ds{-a}{b}\rtimes \sigma_c) 
            +
            1 \cdot L(\ds{-b}{c}\rtimes \sigma_a)
        \\
        \ds{\p}{c}\rtimes \sigma_{a,b}^+
        &\geq
           1 \cdot L(\ds{-a}{b}\rtimes \sigma_c) 
            +
            1 \cdot L(\ds{-b}{c}\rtimes \sigma_a).
    \end{align*}
\end{subproposition}
\begin{proof}
   Both claims follow by Lemmas 8.4 and  9.1 of 
   \cite{ciganovic-glasnik}.
\end{proof}

\section{Composition factors of
    $ \ds{\p}{c} \rtimes \sigma_{a,b}^+$
    and
    $ \ds{\p}{c} \rtimes \sigma_{a,b}^-$.
    }
\label{faktori-malih-reprezentacija}

Here we take in consideration 
Propositions 
\ref{svi-diskretni-u-pola-c--a-b-sigma}
and 
\ref{nediskretni-kandidati-u-pozitivnoj-C}
as well as all corollaries 
of Section 
\ref{confirming-non-tempered subquotients}
to immediately obtain composition factors of above representations.

\begin{theorem}
\label{faktori-od-pola-c-sigma-a-b-pm}
    We have in $R(G)$
\begin{align*}
\ds{\p}{c} \rtimes \sigma_{a,b}^+
=
&L(\ds{-a}{c}\rtimes \sigma_b)
+
L(\ds{-a}{b}\rtimes \sigma_c)
+
L(\ds{-b}{c}\rtimes\sigma_a)
+
\\
&L(\ds{\p}{c} \rtimes \sigma_{a,b}^+)
+
L(\ds{\p}{a}\rtimes \sigma_{b,c}^+) 
+
L(\ds{\p}{b}\rtimes \sigma_{a,c}^+)
+
\\
&\sigma_{b,c,a}^+ + \sigma_{b,c,a}^-,
\\
\\
\ds{\p}{c} \rtimes \sigma_{a,b}^-
=
&L(\ds{\p}{c} \rtimes \sigma_{a,b}^-)
+
L(\ds{\p}{a}\rtimes \sigma_{b,c}^-) 
+
L(\ds{\p}{b}\rtimes \sigma_{a,c}^-)
+
\\
& \sigma_{a,b,c}^-.
\end{align*} 
\end{theorem}

\section{Multiplicity of 
    $L(\ds{-b}{c}\times \ds{\p}{a} \rtimes \sigma$)}
\label{section-multiplicitet-L--b-c-pola-a}

Here we consider multiplicity
of  the above subquotient
in 
$
\ds{\p}{c}
    \times
        \ds{-a}{b}  
                \rtimes \sigma$,
since it appears in two kernels in
the decomposition in Section \ref{dekompozicija1}.

\begin{lemma}
\label{ds(-a,b)otimes-sigma-c-u-LANG}
    We have with maximum multiplicity:
    \begin{align*}
        \mu^*(L(\ds{-b}{c}\times \ds{\p}{a} \rtimes \sigma))
        \geq 
            1 \cdot \ds{-a}{b}\otimes \sigma_c
    \end{align*}
\end{lemma}
\begin{proof}
We can check multiplicity 
one of 
$\ds{-a}{b}\otimes \sigma_c$ in 
    $\mu^*(\ds{\p}{a} \rtimes L(\ds{-b}{c} \rtimes \sigma))$,
using \eqref{jacquet-langlandsov-kvocijent}.
By Corollary 4.1, we have in $R(G)$:
\[
    \ds{\p}{a} \rtimes L(\ds{-b}{c} \rtimes \sigma)
    =
    L(\ds{-b}{c}\rtimes \sigma_a)
        + 
            L(\ds{-b}{c}\times \ds{\p}{a} \rtimes \sigma).
\]

To see that 
$ \ds{-a}{b}\otimes \sigma_c \nleq \mu^*(L(\ds{-b}{c}\rtimes \sigma_a))$, check multiplicity one of 
$\ds{-a}{b}\otimes \sigma_c$
in 
$\ds{-b}{c}\rtimes \sigma_a$ and observe that it comes from $\sigma_{a,b,c}^+$.
\end{proof}

Using 
\eqref{komnozenje} 
and  
\eqref{jacquet-segment},
it is not hard to check the following
\begin{lemma} 
\label{multi-ds(-a,b)otimes-sigma-c}
    We have with maximum multiplicity:
    \begin{align}
    \label{multi-ds(-a,b)otimes-sigma-c-1}
        \mu^*(\ds{-b}{c}\times \ds{\p}{a} \rtimes \sigma))
        \geq 
            2 \cdot \ds{-a}{b}\otimes \sigma_c,
          \\
        \label{multi-ds(-a,b)otimes-sigma-c-2}
          \mu^*(\ds{-a}{c}\times \ds{\p}{b} \rtimes \sigma))
        \geq 
            2 \cdot \ds{-a}{b}\otimes \sigma_c,
        \\
        \label{multi-ds(-a,b)otimes-sigma-c-3}
        \mu^*(\ds{-a}{b}\times \ds{\p}{c} \rtimes \sigma))
        \geq 
            4 \cdot \ds{-a}{b}\otimes \sigma_c,
        \\
        \label{multi-ds(-a,b)otimes-sigma-c-4}
        \mu^*( \ds{\p}{c} \rtimes \sigma_{a,b}^+)
        \geq 
            2 \cdot \ds{-a}{b}\otimes \sigma_c,
    \end{align}
\end{lemma}

\begin{lemma}
\label{ds(-a,b)otimes-sigma-c-u-LANG-2}
We have with maximum multiplicity:
\begin{align*}
\mu^*( 
    L(\ds{\p}{c} \rtimes \sigma_{a,b}^+) )
         \geq 
            1 \cdot \ds{-a}{b}\otimes \sigma_c
\end{align*}
\end{lemma}
\begin{proof}
By Proposition 
\ref{faktori-od-pola-c-sigma-a-b-pm}  and Theorem 
\ref{faktori-pola-b-ds-a-c}
all irreducible subquotients of
\\
$\ds{\p}{c} \rtimes \sigma_{a,b}^+$,
except 
$L(\ds{\p}{c} \rtimes \sigma_{a,b}^+)$
and
$\sigma_{a,b,c}^-$,
are contained in 
$\ds{-a}{c}\times \ds{\p}{b} \rtimes \sigma$.
Now 
\eqref{multi-ds(-a,b)otimes-sigma-c-2}
and
\eqref{multi-ds(-a,b)otimes-sigma-c-4}
imply the claim, since exactly $\sigma_{a,b,c}^+$
and,
by Lemma
\ref{ds(-a,b)otimes-sigma-c-u-LANG},
$L(\ds{-b}{c}\times \ds{\p}{a} \rtimes \sigma)
\nleq 
\ds{\p}{c} \rtimes \sigma_{a,b}^+
$
contribute with
$\ds{-a}{b}\otimes \sigma_c$
in 
$\mu^*(\ds{-a}{c}\times \ds{\p}{b} \rtimes \sigma)$.
\end{proof}

Using \eqref{multi-ds(-a,b)otimes-sigma-c-4} and Lemmas 
\ref{ds(-a,b)otimes-sigma-c-u-LANG}
and
\ref{ds(-a,b)otimes-sigma-c-u-LANG-2}
have
\begin{lemma}
\label{lista-potkvoc-sa-ds(-a,b)otimes-sigma-c}
    Let $\pi$ be an irreducible subquotient of 
    $\ds{\p}{c} \times \ds{-a}{b} \rtimes \sigma$.
    Then
$\mu^*(\pi)\geq 
\ds{-a}{b}\otimes \sigma_c$
if and only if
$\pi$ is one of the following
\begin{align*}
    \sigma_{a,b,c}^+,
    \quad
    \sigma_{a,b,c}^-,
    \quad
    L(\ds{-b}{c}\times \ds{\p}{a} \rtimes \sigma),
    \quad
    L(\ds{\p}{c} \rtimes \sigma_{a,b}^+).    
\end{align*}
\end{lemma}

\begin{proposition}
\label{mult-L5}
    We have in $R(G)$ with maximum multiplicity
    \begin{align*}
   \ds{\p}{c} \times \ds{-a}{b} \rtimes \sigma
   \geq 
   1 \cdot
    L(\ds{-b}{c}\times \ds{\p}{a} \rtimes \sigma)
    \end{align*}
\end{proposition}
\begin{proof}
Follows directly from
\eqref{multi-ds(-a,b)otimes-sigma-c-3}
and
Lemma \ref{lista-potkvoc-sa-ds(-a,b)otimes-sigma-c}

\end{proof}

\section{Composition factors of
$   \ds{\p}{c}\times
    \ds{-a}{b}\rtimes \sigma    $ and the first filtration  }

\label{faktori-velike-reprezentacije}

Finally we determine composition factors of the induced representation.

\begin{theorem}
\label{faktori-pola-c-ds-a-b}
We have in $R(G)$
\begin{align*}
\ds{\p}{c}\times
    \ds{-a}{b}\rtimes \sigma
 =&L(\ds{\p}{a} \times
    \ds{-b}{c}\rtimes \sigma)+
 \\
 &
 L(\ds{\p}{b} \times
 \ds{-a}{c} \rtimes \sigma)+
 \\
 &
 L(\ds{\p}{c}\times
    \ds{-a}{b}\rtimes \sigma)
\\
&
+
\\
&L(\ds{-a}{c}\rtimes\sigma_{b})
    +L(\ds{-a}{b}\rtimes \sigma_c)
        +L(\ds{-b}{c}\rtimes \sigma_a)
\\
&
+
\\
&
L(\ds{\p}{a}\rtimes \sigma_{b,c}^+)
    +
        L(\ds{\p}{b} \rtimes \sigma_{a,c}^+)
            +
                L(\ds{\p}{c}\rtimes \sigma_{a,b}^+)
+
\\
&
L(\ds{\p}{a}\rtimes \sigma_{b,c}^-)
    +
        L(\ds{\p}{b}\rtimes \sigma_{a,c}^-)
            +
                L(\ds{\p}{c}\rtimes \sigma_{a,b}^-)
\\
&
+
\\
&\sigma_{a,b,c}^+ 
+
\sigma_{a,b,c}^-
+
\sigma_{b,c,a}^-.
\end{align*}
\end{theorem}

\begin{proof}
    Theorems
    \ref{muic-A3},
    \ref{ciganovic-A2},
    \ref{faktori-pola-b-ds-a-c}
    and
    \ref{faktori-od-pola-c-sigma-a-b-pm}
    describe kernels $K_i$ in \eqref{long-intertwining},
    and thus determine all composition factors.
    Multiplicity one for discrete series is proved in 
    Propostiton
    \ref{svi-diskretni-u-pola-c--a-b-sigma}.
    For the other subquotients, 
    appearing in two or more kernels $K_i$,
    multiplicity one is proved in Propostions
    \ref{mult-L1}, \ref{mult-L2}, \ref{mult-L3}, \ref{mult-L4}
    and \ref{mult-L5}.
\end{proof}

For the sake of the completeness, we say here some more information about kernels from
Section \ref{dekompozicija1}. Using notation there, let us denote
\[
    k_1=K_1, \quad k_{i}= K_{i}\cap \textrm{Im}(f_{i-1}\circ\cdots \circ f_1), \quad i=2,3,4.
\]
\begin{lemma}
    We have
    \begin{align*}
        k_2\cong \sigma_{a,b,c}^-, \quad 
        k_3\cong \{0\}, \quad
        k_4\cong    L(\ds{\p}{c}\rtimes \sigma_{a,b}^+)+
                    L(\ds{\p}{c}\rtimes \sigma_{a,b}^-).
    \end{align*}
\end{lemma}
\begin{proof} Theorem
\ref{faktori-pola-b-ds-a-c} gives composition series of $K_1$. 
Proposition \ref{druga}        
and Theorems \ref{ciganovic-A2},
\ref{faktori-od-pola-c-sigma-a-b-pm}
and \ref{faktori-pola-c-ds-a-b} give composition factors of 
$K_2$, $K_3$, $K_4$ and $\psi$. 
All factors of $K_2$, except $\sigma_{a,b,c}^-$, are in $K_1$ giving the first equation. All factors in $K_3$ are in either $K_1$ or $K_2$, giving the second equation. Similarly goes for the last equation.
\end{proof}

This gives us a filtration of $\psi$, but in Section \ref{main result} we obtain a more precise result.

\begin{corollary}
There exists a filtration 
$\{V_i\}$
of
 $\ds{\p}{c}\times \ds{-a}{b}\rtimes \sigma$ 
 such that
\begin{align*}
V_1 \cong &
    \sigma_{b,c,a}^+ 
    +L(\ds{\p}{a}\rtimes \sigma_{b,c}^-),
        \\
V_2/V_1 \cong &
    L(\ds{\p}{a}\rtimes \sigma_{b,c}^+)
    +
    L(\ds{-a}{b}\rtimes \sigma_c)
    +
    L(\ds{\p}{b}\rtimes \sigma_{a,c}^-)
        +
        L(\ds{-b}{c}\rtimes \sigma_a),
        \\
V_3/V_2\cong & 
    L(\ds{\p}{b} \rtimes \sigma_{a,c}^+)+
            L(\ds{-a}{c}\rtimes\sigma_{b})+
        \sigma_{b,c,a}^- + 
                            L(\ds{-b}{c}\times
                                \ds{\p}{a} 
                                    \rtimes \sigma),
        \\
V_4/V_3\cong  &
    L(\ds{-a}{c}\times \ds{\p}{b}\rtimes \sigma),    
\\
V_5/V_4\cong  & \sigma_{a,b,c}^-,
\\
V_6/V_5\cong &
    L(\ds{\p}{c}\rtimes \sigma_{a,b}^+)+
                    L(\ds{\p}{c}\rtimes \sigma_{a,b}^-).
\end{align*}

\end{corollary}

\section{Composition series of
    $ \ds{\p}{c} \rtimes \sigma_{a,b}^+$
    and
    $ \ds{\p}{c} \rtimes \sigma_{a,b}^-$.
    }

\label{kompozicijski-niz-malih}

Here we determine composition series of kernel $K_4$ from  Section \ref{dekompozicija1}.

\begin{theorem}
There exists a filtration $\{ V_i\}$ of 
    $\ds{\p}{c} \rtimes \sigma_{a,b}^+$
    where
\begin{align*}
    V_1=&\sigma_{b,c,a}^+ ,
        \\
    V_2/V_1 \cong &L(\ds{\p}{a}\rtimes \sigma_{b,c}^+)
        +
        L(\ds{-a}{b}\rtimes \sigma_c)
        +
        L(\ds{-b}{c}\rtimes \sigma_a),
        \\
     V_3/V_2 \cong &L(\ds{\p}{b} \rtimes \sigma_{a,c}^+)+
            L(\ds{-a}{c}\rtimes\sigma_{b})+
                \sigma_{b,c,a}^- ,
        \\
    V_4 / V_3 \cong &L(\ds{\p}{c} \rtimes \sigma_{a,b}^+).                    
\end{align*}
\end{theorem}

\begin{proof} 
By Theorem 
\ref{faktori-pola-c-ds-a-b}
we have multiplicity one representations:
    \begin{align*}
        \ds{\p}{c} \rtimes \sigma_{a,b}^+
        \hookrightarrow
        \ds{\p}{c}\times
        \ds{-a}{b}\rtimes \sigma     
        \hookleftarrow
        \ds{\p}{b}\times\ds{-a}{c}\rtimes \sigma.
    \end{align*}
By Theorems
\ref{faktori-pola-b-ds-a-c} and
\ref{faktori-od-pola-c-sigma-a-b-pm}
all irreducible subquotents of
$\ds{\p}{c} \rtimes \sigma_{a,b}^+$, except its Langlands quotient,
appear as irreducible subquotients in 
$\ds{\p}{b}\times\ds{-a}{c}\rtimes \sigma$.
We conclude that
the unique maximal submodule $M$ of 
\( \ds{\p}{c} \rtimes \sigma_{a,b}^+\)
is isomorphic to a subrepresentation of 
\( \ds{\p}{b}\times\ds{-a}{c}\rtimes \sigma\).
Further, by filtrations 
of 
\( \ds{\p}{b}\times\ds{-a}{c}\rtimes \sigma\) 
(Theorem 14.1 of \cite{ciganovic-glasnik})
and
\( \ds{\p}{b} \rtimes \sigma_{a,c}^-\)
(Proposition 10.2. of \cite{ciganovic-glasnik})
we have 
\[
    M\hookrightarrow 
        \ds{\p}{b}\times\ds{-a}{c}\rtimes \sigma 
            /
                \ds{\p}{b} \rtimes \sigma_{a,c}^-
\]
and  the claim follows.
\end{proof}

Now we consider composition series of 
$ \ds{\p}{c} \rtimes \sigma_{a,b}^-$. First we need a lemma.

\begin{lemma} We have an embedding
\label{ulaganje pola b - u pola c -}
\begin{equation} 
    \ds{\p}{b} \rtimes \sigma_{a,c}^-
        \hookrightarrow
            \ds{\p}{c} \rtimes \sigma_{a,b}^-.
\end{equation}    
\end{lemma}
\begin{proof}
By Theorem \ref{faktori-pola-c-ds-a-b} 
we have multiplicity one representations:
    \begin{align*}
        \ds{\p}{b} \rtimes \sigma_{a,c}^-
            &\hookrightarrow
                \ds{\p}{b} \times \ds{-a}{c} \rtimes \sigma
                    \\
                    &\hookrightarrow
                        \ds{\p}{c} \times \ds{-a}{b} \rtimes \sigma
                            \hookleftarrow
                                \ds{\p}{c} \rtimes \sigma_{a,b}^-.
    \end{align*}
    By Lemma \ref{faktori-pola-b-sigma -a c pm}
    and Theorem \ref{faktori-od-pola-c-sigma-a-b-pm},
    all irreducible subquotients of 
    $\ds{\p}{b} \rtimes \sigma_{a,c}^-$ 
    appear in composition factors of
    $\ds{\p}{c} \rtimes \sigma_{a,b}^-$,
    and the claim follows.
\end{proof}
By Lemmas \ref{faktori-pola-b-sigma -a c pm}
and \ref{ulaganje pola b - u pola c -}
we have composition series
\begin{theorem}
    There exists a filtration $\{ V_i\}$ of 
    $\ds{\p}{c} \rtimes \sigma_{a,b}^-$
    where
    \begin{alignat*}{3}
        &V_1 \cong        L(\ds{\p}{a}\rtimes \sigma_{b,c}^-),
        \quad \quad
        &V_2/V_1 &\cong   L(\ds{\p}{b} \rtimes \sigma_{a,c}^-),
        \\
        &V_3/V_2 \cong         \sigma_{a,b,c}^-,
        \quad \quad
        &V_4/V_3 &\cong   L(\ds{\p}{c} \rtimes \sigma_{a,b}^-).
    \end{alignat*}
\end{theorem}

\section{Composition series of 
$\ds{\p}{c}\times \ds{-a}{b} \rtimes \sigma$}   
\label{main result}

Now we have the main result.

\begin{theorem} 
Let $\psi=\ds{\p}{c}\times \ds{-a}{b} \rtimes \sigma$ and define representations
\begin{align*}
    W_1=&\sigma_{b,c,a}^+ +L(\ds{\p}{a}\rtimes \sigma_{b,c}^-),
        \\
    W_2=&L(\ds{\p}{a}\rtimes \sigma_{b,c}^+)
        +
        L(\ds{\p}{b}\rtimes \sigma_{a,c}^-)
        +
        L(\ds{-b}{c}\rtimes \sigma_a)
        +
        L(\ds{-a}{b}\rtimes \sigma_c),
        \\
     W_3=&
            \sigma_{b,c,a}^- +  \sigma_{a,b,c}^-
            +
            L(\ds{\p}{b} \rtimes \sigma_{a,c}^+)
            +
            L(\ds{-a}{c}\rtimes\sigma_{b})+
             L(\ds{-b}{c}\times \ds{\p}{a} 
                                    \rtimes \sigma),
        \\
      W_4=&L(\ds{\p}{b} \times \ds{-a}{c} \rtimes \sigma)
        + L(\ds{\p}{c} \rtimes \sigma_{a,b}^+)
        +L(\ds{\p}{c} \rtimes \sigma_{a,b}^-),
        \\
    W_5=&L(\psi).
\end{align*}
Then there exists a sequence 
$\{0\}=V_0\subseteq V_1 \subseteq V_2 \subseteq V_3
            \subseteq V_4 \subseteq V_5=\psi$,
such that
\begin{equation*}
V_i/V_{i-1}\cong W_i,\quad  i=1,\ldots,5.
\end{equation*}    
\end{theorem}
\begin{proof}
    Since 
    \[\ds{\p}{b} \times \ds{-a}{c} \rtimes \sigma \hookrightarrow \psi,\] 
    its filtration, 
    Theorem \ref{faktori-pola-b-ds-a-c}, 
    imply existence of $V_1$ and $V_2$.
    Additionaly, composition factors of 
    $\ds{\p}{c} \rtimes \sigma_{a,b}^-$, 
    Theorem \ref{faktori-od-pola-c-sigma-a-b-pm},
    and
        $\ds{\p}{c} \rtimes \sigma_{a,b}^- \hookrightarrow
            \psi$
    imply 
    \[
        \sigma_{a,b,c}^- \hookrightarrow \psi /V_2
                \hookleftarrow 
                    \ds{\p}{b} \times \ds{-a}{c} \rtimes \sigma/V_2.
    \]
    This shows $W_3 \hookrightarrow \psi/V_2$ proving existence of $V_3$.
    Finally, Theorems \ref{faktori-pola-b-ds-a-c}, 
    \ref{faktori-od-pola-c-sigma-a-b-pm}
    and
    \[
            \ds{\p}{c} \rtimes \sigma_{a,b}^\pm \hookrightarrow \psi 
                \hookleftarrow 
                    \ds{\p}{b} \times \ds{-a}{c} \rtimes \sigma.
    \]
    show 
    $W_4 \hookrightarrow \psi/V_3$, proving existence of $V_4$.
\end{proof}



\newpage

\def\cprime{$'$} \def\cprime{$'$}
\providecommand{\bysame}{\leavevmode\hbox to3em{\hrulefill}\thinspace}
\providecommand{\MR}{\relax\ifhmode\unskip\space\fi MR }
\providecommand{\MRhref}[2]{
  \href{http://www.ams.org/mathscinet-getitem?mr=#1}{#2}
}


\begin{thebibliography}{10}



\bibitem{bernstein-zelevinsky:ind-repns-I}
I.~N. Bernstein and A.~V. Zelevinsky, \emph{Induced representations of
  reductive {$p$}-adic groups. {I}}, Ann. Sci. \'{E}cole Norm. Sup. (4)
  10 (1977), no.~4, 441-472.






\bibitem{bbosnjak:ladder}
B. Bošnjak,
\emph{Representations induced from cuspidal and ladder representations 
of classical $p$-adic groups.}
Proc. Amer. Math. Soc. 149 (2021), no. 12, 5081-5091.





\bibitem{ciganovic}
I.~Ciganovi\'c, 
\emph{Composition series of a class of induced representations, a case of one half cuspidal reducibility}, Pacific J. Math. 296 (2018),  no. 1, 21-30.


\bibitem{ciganovic-discrete}
I.~Ciganovi\'c,
\emph{Composition series of a class of induced representations built on discrete series},
Manuscripta Math. 170 (2023), no. 1-2, 1–18.


\bibitem{ciganovic-glasnik}
I.~Ciganovi\'c,
\emph{Parabolic induction from two segments, linked under contragredient, with a one half cuspidal reducibility, a special case},
Glas. Mat. Ser. III 59(79)(2024), no. 1, 77–105.


\bibitem{goldberg-tempered}
D.~Goldberg,
\emph{Reducibility of induced representations for $Sp(2n)$ and $SO(2n)$},
Amer J Math, 1994, 116: 1101-1151.


\bibitem{muic-hanzer:langzel}
M. Hanzer, G. Muić, 
\emph{On an algebraic approach to the Zelevinsky
classification for classical $p$-adic groups},
J.Algebra 320 (2008), no. 8, 3206-3231.










\bibitem{jantzen:support}
C. Jantzen, 
\emph{On supports of induced representations for symplectic and 
odd-orthogonal groups}, 
Amer. J. Math., 119 (1997), no. 6, 1213-1262.





\bibitem{degenerate}
Y. Kim, B. Liu, and I. Matić,
\emph{Degenerate principal series for classical and odd GSpin
groups in the general case}, 
Represent. Theory 24 (2020), 403-434.



\bibitem{faktor}
I.~Mati\'c, \emph{On discrete series subrepresentations 
of the generalized principal series}, Glas. Mat. Ser.III 51(71)(2016),
no. 1, 125-152




\bibitem{matic-forum}
I.~Mati\'c
\emph{Representations induced from the Zelevinsky segment and 
discrete series in the half-integral case},
Forum Math. 33 (2021), no. 1, 193-212.




\bibitem{segment}
I.~Mati\'c, M.~Tadi\'c, \emph{On Jacquet modules of representations of segment type}, Manuscripta Math. 147 (2015), no. 3-4, 437-476.


\bibitem{segment-ispravljeno} 
I.~Mati\'c, M.~Tadi\'c, \emph{On Jacquet modules of representations of segment type}, 
\url{http://www.hazu.hr/~tadic/53-JM-segment.pdf}





\bibitem{moeglin}
C.~M{\oe}glin, \emph{Sur la classification des s\'{e}ries discr\`{e}tes des groupes classiques $p$-adiques: param\`{e}tres de Langlands et exhaustivit\'{e}. 
} 
J. Eur. Math. Soc. 
 4 (2002), no. 2, 143-200.


\bibitem{M-V-W}
C. M{\oe}glin, M.-F. Vigneras, and J.-L. Waldspurger, \emph{Correspondence de Howe sur un corps p–adique},
Lecture Notes in Math. 1291, 1987.


\bibitem{tadic-diskretne}
C.~M{\oe}glin, M.~Tadi\'{c}, 
\emph{Construction of discrete series for classical $p$-adic groups,}
J. Amer. Math. Soc. 15 (2002), no. 3, 715-786. 





\bibitem{muic:composition series}
G.~Mui{\'c}, \emph{Composition series of generalized principal series; the case of strongly positive discrete series}, Israel J. Math. 140 (2004), 157-202. 




\bibitem{muic:reducibility principal}
\bysame, \emph{Reducibility of generalized principal series}, Canad. J. Math. 57 (2005), no. 3, 616-647.








\bibitem{tadic:square-integr-correspond-segments}
M.~Tadi{\'c}, \emph{Square integrable representations of classical p-adic groups corresponding to segments}. Represent Theory 3 (1999), 58-89.









\bibitem{tadic:red-par-ind}
\bysame, \emph{On reducibility of parabolic induction}, Israel J. Math.
  \textbf{107} (1998), 29-91.


\bibitem{tadic:unitary}
\bysame, \emph{An external approach to unitary 
representations}, Bull. Amer. Math. Soc. (N.S.),
28 (1993), no. 2, 215-252.



\bibitem{tadic:regular-square}
\bysame, \emph{On regular square integrable representations of $p$-adic groups}, Amer. J. Math. 120 (1998), no. 1, 159-210.



\bibitem{tadic-structure}
\bysame, \emph{Structure arising from induction and Jacquet modules of representations of classical $p$-adic groups}, J. Algebra 177 (1995), no. 1, 1-33.


\bibitem{tadic:family}
\bysame,
\emph{A family of square integrable representations of classical $p$-adic groups in the case of general half-integral reducibilities} Glas. Mat. Ser. III 37(57) (2002), no. 1, 21-57.



\bibitem{tadic:tempered}
\bysame,
\emph{On tempered and square integrable representations of classical $p$-adic groups},
Sci. China Math., 56 (2013), pp. 2273–2313.


\bibitem{tadic:corank}
\bysame,\emph{ Unitarizability in corank three for classical p-adic groups},
 Mem. Amer. Math.
Soc., 286 (2023), pp. vii+120.



\bibitem{zelevinsky:ind-repns-II}
A.~V. Zelevinsky, \emph{Induced representations of reductive $p$-adic groups.
  {II}. On irreducible representations of $GL(n)$}, Ann. Sci. \'{E}cole
  Norm. Sup. (4) \textbf{13} (1980), no.~2, 165-210.



\end{thebibliography}
\end{document}